\documentclass[10pt]{amsart}
\usepackage{amsmath,mathrsfs}
\usepackage{amssymb}
\usepackage{latexsym}
\usepackage{amsfonts}
\usepackage{graphicx}
\usepackage{psfrag}

\newtheorem{Pa}{Paper}[section]
\newtheorem{theorem}[Pa]{{\bf Theorem}}
\newtheorem{lemma}[Pa]{{\bf Lemma}}

\newtheorem{proposition}[Pa]{{\bf Proposition}}
\newtheorem{observation}[Pa]{{\bf Observation}}

\newtheorem{example}[Pa]{{\bf Example}}

\def\C{\mathbb C}
\def\R{{\mathbb R}}

\def\@tempb{saamsart}
\usepackage{hyperref}
\begin{document}
\bibliographystyle{plain}
\thispagestyle{empty}
\author[D. Alpay]{Daniel Alpay}
\address{(DA) Department of Mathematics
\newline
Ben Gurion University of the Negev \newline P.O.B. 653,
\newline
Be'er Sheva 84105, \newline ISRAEL} \email{dany@math.bgu.ac.il}
\author[I. Lewkowicz]{Izchak Lewkowicz}
\address{(IL) Department of Electrical Engineering
\newline
Ben Gurion University of the Negev \newline P.O.B. 653,
\newline
Be'er Sheva 84105, \newline ISRAEL }
\email{izchak@ee.bgu.ac.il}
\date{}
\title[
minimal realization of families of systems]
{minimal state space realization, static output feedback
and matrix completion of families of systems}
\thanks{
This research was supported in part by the Bi-national Science
Foundation grant no. 2010117.
D. Alpay thanks the
Earl Katz family for endowing the chair
which supported his research.}

\def\squarebox#1{\hbox to #1{\hfill\vbox to #1{\vfill}}}

\noindent

\noindent

\begin{abstract}\quad
We here characterize the minimality of realization of arbitrary
linear time-invariant dynamical systems through (i) intersection
of the spectra of the realization matrix and of the
corresponding state submatrix and (ii) moving the
poles by applying static output feedback.

\noindent
In passing, we introduce, for a given square matrix $A$, a
parameterization of all matrices $B$ for which the pairs
$A, B$ are controllable. In particular, the minimal rank
of such $B$ turns to be equal to the smallest geometric
multiplicity among the eigenvalues of $A$.

\noindent
Finally, we show that the use of a (not necessarily square)
realization matrix $L$ to examine minimality of realization,
is equivalent to the study of $L_{\rm sq}$, a smaller
dimensions, square realization matrix, which in turn is
linked to realization matrices obtained as polynomials in
$L_{\rm sq}$. Namely a whole family of systems.
\end{abstract}
\maketitle

\section{Introduction}
\label{sec:Introduction}
\setcounter{equation}{0}

We first recall in the concept of state space realization. Let
$F(s)$ be a $p\times m$-valued rational function ($m$ inputs
and $p$ outputs in engineering terminology), analytic at
infinity, i.e. $\lim\limits_{s~\rightarrow~\infty}~F(s)~$
exists. Then, $F(s)$ admits a state space realization
\begin{equation}\label{StateSpace}
\begin{matrix}
F(s)=C(sI-A)^{-1}B+D&~&~&
L:={\footnotesize
\left(\begin{array}{r|r}A&~B\\ \hline C&~D\end{array}
\right)
}
\end{matrix}
\end{equation}
with $~A\in\C^{n\times n}$, $~B\in\C^{n\times m}$,
$~C\in\C^{p\times n}~$ and $~D\in\C^{p\times m}$, namely,
\mbox{$L\in\C^{(n+p)\times(n+m)}$}.
\vskip 0.2cm

For convenience, throughout this work the formulation is for
complex matrices. However, up to complex conjugation of
spectrum, all results hold for real matrices.
\vskip 0.2cm

A realization is called~ {\em minimal} whenever $~n$, the
dimension of $A$ in \eqref{StateSpace}, is the smallest possible,
see e.g. \cite[Definition 6.5.9]{So}. It is
then equal to the McMillan degree of $F(s)$, see e.g.
\cite[6.5-9]{Ka}, \cite[Remark 6.7.4]{So}.
\vskip 0.2cm

The issue of minimality of realization is fundamental,
see e.g. \cite[Sections 2.5, 2.6.2]{DP}, 
\cite[Section 2.4, Theorem 6.2-3]{Ka},
\cite[Sections 6.2, 6.3]{LR}. Typically,
one is first interested in the question of whether or not a given
realization is minimal and if not, to find ways to extract out of
the given realization, a minimal one. For a survey of
works addressing the second question see \cite{DS}. We here focus
on the first problem.
\vskip 0.2cm

For a (possibly rectangular) system $F(s)$ admitting state space
realization \eqref{StateSpace}, we shall also find it convenient
to consider the~ {\em squared} realization matrix $L_{\rm sq}$
without altering $A$ and thus preserving $n$. The naive
(=inflating) version
is obtained as follows: If $m>p$ by simply adding rows of zeros
to the bottom $L$ until it is \mbox{$(n+m)\times(n+m)$}, and if
$p>m$ by adding columns of zeros to the right side of $L$ until
it is \mbox{$(n+p)\times(n+p)$}.
This approach preserves the associated controllable and observable
subspaces. In particular, if exists, minimality of realization is
maintained. As mentioned, here the dimension of $L_{\rm sq}$ is
equal to the larger \mbox{dimension of the original $L$.}
\vskip 0.2cm

In Lemma \ref{L:TruncSq} below, we show that if a realization
$L$ in \eqref{StateSpace} is minimal, one can multiply the
corresponding rational function $F(s)$ by constant matrices:
$T_c$ from the left and $T_b$ from the right, so that the
resulting $T_cF(s)T_b$ is of the smallest possible dimensions
(which turn to be square) while in the resulting realization,
$L_{\rm sq}$, both the $A$ matrix and minimality, are preserved.
\vskip 0.2cm

Recall that for a given matrix $A\in\C^{n\times n}$, the~
{\em geometric multiplicity}~ of \mbox{$\lambda\in{\rm spect}(A)$}
is the number of linearly independent eigenvectors associated
with this $\lambda$, see e.g. \cite[Definition 1.4.3]{HJ1}. 
The largest geometric multiplicity among the eigenvalues of
$A$ will be denoted by $\alpha(A)$.
\vskip 0.2cm

The three main results of this work are Proposition \ref{L:TruncSq},
Theorem \ref{Th:Main} and Proposition \ref{RationaL} are stated in
the Introduction.

\begin{proposition}\label{L:TruncSq}
Denote by $\alpha=\alpha(A)$ the largest geometric multiplicity among
the eigenvalues of a given a matrix $A\in\C^{n\times n}$.
\vskip 0.2cm

The following is true.

\begin{itemize}
\item[(i)~~~]{} If for a given a matrix $B\in\C^{n\times m}$ the pair
$A, B$ is controllable, then
\[
{\rm rank}(B)\geq\alpha.
\]
Moreover, there exists a full
rank matrix ~\mbox{$T_b\in\C^{m\times\alpha}$}
(and \mbox{$\hat{B}:=BT_b\in\C^{n\times\alpha}$} is
of a full rank)
so that the pair $A, \hat{B}$ is controllable.
\vskip 0.2cm

\item[(ii)~~]{} If for a given a matrix $C\in\C^{p\times n}$ the pair
$A, C$ is observable, then
\[
{\rm rank}(C)\geq\alpha.
\]
Moreover, there exists a full rank matrix~
\mbox{$T_c\in\C^{\alpha\times p}$} (and
\mbox{$\hat{C}:=T_cC\in\C^{\alpha\times n}$} is of a full
rank) so that the pair $A, \hat{C}$ is observable.
\vskip 0.2cm

\item[(iii)~]{} Given a $~p\times m$-valued rational function $F(s)$
and the corresponding \mbox{$(n+p)\times(n+m)$} realization $L$ as in 
\eqref{StateSpace}. If the realization $L$ is minimal,
there exist full rank matrices~
\mbox{$T_b\in\C^{m\times\alpha}$} and~
\mbox{$T_c\in\C^{\alpha\times p}$} so that
the $(n+\alpha)\times(n+\alpha)$ matrix
\begin{equation}\label{eq:sq}
L_{\rm sq}:=\left(\begin{smallmatrix}I_n&0\\ 0&T_c
\end{smallmatrix}\right)L\left(\begin{smallmatrix}
I_n&0\\ 0&T_b\end{smallmatrix}\right)=
{\footnotesize \left(\begin{array}{c|c}A&~\hat{B}\\ 
\hline \hat{C}&~T_cDT_b\end{array}\right)}
\end{equation}
is a minimal realization of the ~$\alpha\times\alpha$-valued
rational function
\end{itemize}
\[
F_{\rm sq}(s):=T_c{F}(s)T_b~.
\]
\end{proposition}

To proceed, we need further notation.
For a system $F(s)$ admitting state space realization
\eqref{StateSpace}, we shall find it convenient to consider
the~ {\em associated} system, $\hat{F}(s)$, with zero at infinity
(``strictly proper" in engineering jargon) i.e.
\begin{equation}\label{AssociatedSystem}
\hat{F}(s):=F(s)-F(\infty)=C(sI-A)^{-1}B.
\end{equation}
Thus, the corresponding
$(n+p)\times(n+m)$ realization matrix is of the form
\begin{equation}\label{eq:HatL}
\hat{L}={\footnotesize\left(\begin{array}{r|r}
A&B\\ \hline C&0\end{array}\right)}.
\end{equation}
Note that in both realizations, $L$ in \eqref{StateSpace}
and $\hat{L}$ in \eqref{eq:HatL} the associated controllable
(observable) subspace is identical. In particular, minimality
of these realizations is equivalent.
\vskip 0.2cm

Recall that applying a static output feedback, see
e.g. \cite[Section 3.1]{Ka}, to an
input-output (associated) system \mbox{$y(s)=\hat{F}(s)u(s)$}
means taking $u=Ky+u'$ with $u'$ an auxiliary input and
a (constant) $m\times p$ matrix $K$. The resulting
closed loop system is \mbox{$y(s)=F_{\rm c.l.}(s)u'(s)$}
with \mbox{$F_{\rm c.l.}=(I_p-\hat{F}K)^{-1}\hat{F}$}.
The corresponding closed loop realization matrix $L_{\rm cl}$ is
\begin{equation}\label{Lcl}
L_{\rm c.l.}=
{\footnotesize
\left(\begin
{array}{c|c}A_{\rm cl}&~~B\\ \hline
C&~~0\end{array}\right)}\quad\quad\quad
A_{\rm cl}:=A+BKC.
\end{equation}
\vskip 0.2cm

We can now state our second main result which is a novel
characterization of a minimal realization.

\begin{theorem}\label{Th:Main}
Let $L\in\C^{(n+p)\times (n+m)}$ be a realization of a
$p\times m$-valued rational function $F(s)$,
see \eqref{StateSpace}.

Let $\hat{L}\in\C^{(n+p)\times (n+m)}$ be a realization of
the $p\times m$-valued associated system $\hat{F}(s)$
see \eqref{AssociatedSystem}.
A static output feedback $u=Ky+u'$ (with $K$ $m\times p$
constant) is applied to $\hat{F}(s)$ so that the
closed loop system matrix is as in \eqref{Lcl}.

Let also $L_{\rm sq}$ be a realization of the corresponding
squared system $F_{\rm sq}(s)$.

For $A$ in \eqref{StateSpace} or in \eqref{eq:HatL}, for
some $~q\in[1, n]$ denote
\[
{\rm spect}(A)=\{\lambda_1~,~\ldots~,~\lambda_q\}
\]
with $\lambda_j$ distinct.
\vskip 0.2cm

\noindent
The following are equivalent.
\vskip 0.2cm

\begin{itemize}

\item[(i)~~~]{} The realization of $F(s)$ is minimal.

\item[(ii)~~]{} There exists a (constant) $K\in\C^{m\times p}$
so that in \eqref{Lcl},
\begin{equation}\label{eq:AclA}
{\rm spect}(A)\bigcap{\rm spect}(A_{\rm cl})=\emptyset .
\end{equation}

\item[(iii)~]{} 
Consider the realization $L$ of $F(s)$ in \eqref{StateSpace} with
$m=p$ and assume that $B$ and $C$ are of a full rank. (If this
was not the case take the reduced dimension squared counterpart
$L_{\rm sq}$ and $F_{\rm sq}(s)$ in \eqref{eq:sq}).\quad
Then,
\begin{equation}\label{L_j}
\bigcap\limits_{j=1}^q{\rm spect}
\left(\begin{smallmatrix}A&B\\ C&(\lambda_j-\epsilon)I_p\end{smallmatrix}\right)
\bigcap{\rm spect}(A)=\emptyset ,
\end{equation}
for $\lambda_j\in{\rm spect}(A)$ and 
some sufficiently small $\epsilon>0$.

\item[(iv)~]{}
Consider the realization $L$ of $F(s)$
in \eqref{StateSpace} with $m=p$ (if $m\not=p$ take a squared
counterpart $L_{\rm sq}$ and $F_{\rm sq}(s)$).
\quad
Then,
\begin{equation}\label{NotSpectASpectL}
\bigcap\limits_{D\in\C^{p\times p}}{\rm spect}
\left(\begin{smallmatrix}
A&B\\ C&D\end{smallmatrix}\right)
\bigcap{\rm spect}(A)=\emptyset .
\end{equation}
\end{itemize}
\end{theorem}
\vskip 0.2cm

The above result says that one can study minimality of realization
by examining the spectrum of the matrix $L$ and of its submatrix $A$.
Now, in linear algebra it is natural to discuss families of
matrices sharing the same properties. It is less common to address
families of systems. However, the description of a system by the
realization matrix $L$ in \eqref{StateSpace}, suggests us to
explore ~{\em families}~ of minimal/non-minimal realizations.
\vskip 0.2cm

Before introducing the result, recall in the following. Let $\phi(s)$
be an arbitrary (possibly complex) scalar rational function
(possibly with pole at infinity or even a polynomial). Let
$L$ be an arbitrary $k\times k$ matrix whose eigenvalues differ
from the poles of $\phi(s)$. Then the  $k\times k$ matrix $\phi(L)$
is well defined. In fact, from the Cayley-Hamilton Theorem it
follows that there always exists a polynomial $\psi(s)$, of
degree of at most $k-1$, so that $\phi(L)=\psi(L)$, see e.g.
\cite[Section 2.4]{HJ1}. Thus, we address only polynomials.
\vskip 0.2cm

Although typically one is interested in $F(s)$ (realized by $L$),
but not in the rational function realized by $\psi(L)$, the
result below asserts their interdependence.

\begin{proposition}\label{RationaL}
Let $L\in\C^{(n+m)\times (n+m)}$ be a realization of a
$m\times m$-valued rational function $F(s)$, see \eqref{StateSpace}.

For an arbitrary (scalar) polynomial $\psi(s)$ consider the
$(n+m)\times (n+m)$ matrix $\psi(L)$ as a realization matrix,
i.e.
\mbox{$\psi(L):={\footnotesize\left(\begin{array}{c|c}
\tilde{A}&\tilde{B}\\ \hline
\tilde{C}&\tilde{D}\end{array}\right)}$}.

\begin{itemize}
\item[(a)~~]{}
Consider the following statements

(i) ${\rm spect}(\tilde{A})\bigcap{\rm spect}(\psi(L))=\emptyset$.

(ii) The realization $ \psi(L)$ is minimal.

(iii) The realization $L$ is minimal.

\noindent
Then (i)~ implies~ (ii) which in turn implies~ (iii).
\vskip 0.2cm

\item[(b)~~]{}If the realization $L$ is not minimal, the realization
$\psi(L)$ is not minimal.
\end{itemize}
\end{proposition}

The outline of the paper is as follows. In Section
\ref{sec:Motivation} we motivate these results and
in Sections \ref{sec:Background},
\ref{sec:Proof} and \ref{sec:Families} we prove Proposition \ref{L:TruncSq},
Theorem \ref{Th:Main} and Proposition \ref{RationaL} respectively.

\section{Motivation}
\label{sec:Motivation}
\setcounter{equation}{0}

In this section we put Theorem \ref{Th:Main} in a broader
perspective.
\vskip 0.2cm

{\bf 1.~ Output feedback}\quad Although typically
stated differently, the following is well known,
see e.g. \cite[Section 2.4.3]{DP}, \cite[Section 4.2]{Ka},
\cite[Chapter 7]{So}.

\begin{observation}
Given an associated system $\hat{F}(s)$ and its realization
$\hat{L}$ see \eqref{AssociatedSystem} and \eqref{eq:HatL}.

The loop is closed by applying a state feedback gain to
a (Luenberger) observed state.
\vskip 0.2cm

The realization $\hat{L}$ is minimal if and only if
the closed loop poles may be located anywhere in the
complex plane.
\end{observation}

In engineering jargon, minimal realization enables one to 
place the poles of a closed loop system anywhere in the
complex plane through a ~{\em dynamic}~ output feedback.
\vskip 0.2cm

The simplicity of ~{\em static} output feedback has made it
very attractive. However, exploring its properties turned
out to be challenging, see e.g. \cite{ElGOAiR}, \cite{HL}, 
\mbox{\cite[Section 3.1]{Ka}} and \cite{SADG}.
\vskip 0.2cm

Condition \eqref{eq:AclA} in Theorem \ref{Th:Main} may be viewed
as establishing a precise connection between minimal realization
and ~{\em static} output feedback.
\vskip 0.2cm

{\bf 2}.~ {\bf Static output feedback and realization of the
inverse rational function}\quad Consider a rational function
$F(s)$ and its realization as in \eqref{StateSpace} with
$D=I$, namely,
\mbox{$L={\footnotesize\left(\begin{array}{c|c}A&B\\ \hline
C&I\end{array}\right)}$}.
\vskip 0.2cm

On the one hand, it is straightforward to verify that a
realization of $F^{-1}(s)$, the inverse rational function, 
is given by
\mbox{$L_{\rm inv}={\footnotesize\left(\begin{array}{c|c}
A-BC&B\\ \hline -C&I \end{array}\right)}$}. See e.g. 
\cite[Exer. 2.8(b)]{DP}, \cite[Exer. 2.2-20]{Ka}.
(Obviously, this has nothing to do with $L^{-1}$, the
inverse of the realization matrix,  addressed in Proposition
\ref{RationaL} and Section \ref{sec:Families}).
\vskip 0.2cm

On the one hand, consider now the associated system
\eqref{AssociatedSystem}
and its realization \eqref{eq:HatL}. Applying to it a static output
feedback with $K=-I$ yields in \eqref{Lcl}
\mbox{${L_{\rm cl}}_{|_{K=-I}}={\footnotesize\left(\begin{array}{c|c}
A-BC&B\\ \hline C&0\end{array}\right)}$}.
\vskip 0.2cm

In each of the three systems $L$, $L_{\rm inv}$ and $L_{\rm cl}$,
the associated controllable (or observable) subspace, is identical.
In particular, minimality of the three realization is equivalent.
Thus, Theorem \ref{Th:Main} addresses also the realization of the
inverse system (whenever exists). For example in
\cite[Sections 5, 6]{MG} the authors in fact studied $L$, $L_{\rm inv}$ 
in the context of Linear Fraction Transformation.
\vskip 0.2cm

{\bf 3.~ Matrix completion}\quad
Matrix completion (a.k.a extension) problems have
been of interest in the past 60 years. Many of them can be
casted in the following framework: A part of a matrix is
prescribed, can one complete the missing part so that
the full matrix will poses certain properties, typically
spectral. For a nice survey, see \cite{Cr}. The case where
the upper triangular part is prescribed was addressed in
\cite{BGRS} (and not cited in \cite{Cr}). 
\vskip 0.2cm

As a special case, assume that $L=\left(\begin{smallmatrix}A&B\\
C&*\end{smallmatrix}\right)$, where $A, B, C$ are prescribed
and $*$ stands for a square unprescribed part. Parameterizing
all possible characteristic polynomials of $L$ is known to
be difficult, see comment following \cite[Theorem 46]{Cr}.
Condition \eqref{NotSpectASpectL} in Theorem \ref{Th:Main} can
be seen as answering a more modest question: Under what conditions
can one complete $\left(\begin{smallmatrix}
A&B\\ C&*\end{smallmatrix}\right)$ with $D$ so that the spectra of
the resulting $L$ and $A$ will not (or will always) intersect.
\vskip 0.2cm

{\bf 4.} ~{\bf Pole placement and matrix completion}\quad
The fact that problems matrix completion and
pole placement through feedback, are linked is well
known. See e.g. \cite{KRW} or the Introduction of \cite{RW}. 
The equivalence of \eqref{Lcl} \eqref{NotSpectASpectL}
in Theorem \ref{Th:Main}, falls into this category.
\vskip 0.2cm

{\bf 5}.~ {\bf A spectral PBH test for minimality}\quad
Consider a rational function $F(s)$ and its realization
as in \eqref{StateSpace} or \eqref{AssociatedSystem}. As already
mentioned, the issue of minimality of realization is fundamental.
\vskip 0.2cm

Adopting Kailath's terminology, for given $A\in\C^{n\times n}$,
$B\in\C^{n\times m}$ and $C\in\C^{p\times n}$ ~{\em The
Popov-Belevitch-Hautus (PBH) Rank Tests}, \cite[Theorem 2.16]{DP},
\cite[Theorem 6.2-6]{Ka} say 
that:
A pair $A, B$ is controllable, if and only if
\begin{equation}\label{eq:PBHcont}
{\rm rank}~(\lambda{I}-A\quad B)=n
\quad\quad\quad\quad\forall\lambda\in\C.
\end{equation}
A pair $A, C$ is observable, if and only if
\begin{equation}\label{eq:PBHobs}
{\rm rank}~\left(\begin{smallmatrix}
\lambda{I}-A\\ C\end{smallmatrix}\right)=n
\quad\quad\quad\quad\forall\lambda\in\C.
\end{equation}
It is clear that in \eqref{eq:PBHcont} and \eqref{eq:PBHobs},
without loss of generality, one can confine the search to
\[
\lambda\in{\rm spect}(A).
\]
We here examine two adaptations of these tests:
\vskip 0.2cm

(i) {\em To minimality of realization (without independently
testing for controllability}

\hfill{\em and for observability)}.

(ii) {\em To consider the spectrum of $L_{\rm sq}$ a square
realization matrix.}
\vskip 0.2cm

To this end recall that a realization is minimal if and
only if it is both controllable and observable, see e.g.
\cite[Theorem 2.33]{DP}, \cite[Theorem 6.2-3]{Ka}
\cite[Definition 6.5.3, Theorem 27]{So}.
\vskip 0.2cm

From a combination of the PBH Rank Tests in \eqref{eq:PBHcont}
and \eqref{eq:PBHobs} it follows says that minimal realization
implies that
\begin{equation}\label{eq:Rank}
\min\limits_{\lambda\in\C}~{\rm rank}~\left(\begin{smallmatrix}
\lambda{I}-A&~B\\ C&~0\end{smallmatrix}\right)=n+\min\left(
{\rm rank}(B),~{\rm rank}(C)\right).
\end{equation}
The following example illustrates the fact that formulating the
converse to \eqref{eq:Rank} is more delicate.
Consider the $2\times 1$ rational function of McMillan degree two
\[
F(s)=\left(\begin{matrix}\frac{2\beta\gamma{s}}{s^2-{\alpha}^2}\\~\\
\frac{\gamma(s+2\alpha)}{s+\alpha}\end{matrix}\right)
\quad\quad\quad\quad 0\not=\alpha, \beta, \gamma\in\R.
\]
Its realization is
\[
L=
{\footnotesize
\left(\begin{array}{rr|r}
\alpha&~~0&\beta\\
0&-\alpha&\gamma\\
\hline
\gamma&~~\beta&0\\
0&~~\alpha&\gamma\end{array}\right)}.
\]
Namely, $n=2$, $m=1$ and $p=2$. Although this realization is
minimal, \eqref{eq:Rank} does not hold:
\[
\min\limits_{\lambda\in\C}~{\rm rank}~\left(\begin{smallmatrix}
\lambda{I}-A&~B\\ C&~0\end{smallmatrix}\right)=
{\rm rank}~\left(\begin{smallmatrix}
\lambda{I}-A&~B\\ C&~0\end{smallmatrix}\right)_{|_{\lambda=0}}
=
{\rm rank}~\left(\begin{smallmatrix}
-\alpha&~~0&~~\beta\\~~
0&~~\alpha&~~\gamma\\~~\gamma& \beta&~~0\\~~
0&~~\alpha&~~\gamma\end{smallmatrix}\right)=2.
\]
(indeed, the non-zero vector
$\left(\begin{smallmatrix}
-\beta\\~~\gamma\\-\alpha\end{smallmatrix}\right)$
is in the nullspace of the rightmost matrix).

On the other hand
\[
n+\min\left({\rm rank}(B),~{\rm rank}(C)\right)=2+1=3.
\]
One may view \eqref{NotSpectASpectL} in Theorem \ref{Th:Main} as a
correct extension of the PBH Rank tests to the realization matrix
$L$.
\vskip 0.2cm

{\bf 6.} {\bf Rational functions: Scalar (SISO) versus
matrix-valued (MIMO)
}\quad

\begin{proposition}
Let $L\in\C^{(n+p)\times (n+m)}$ be a realization of a
$p\times m$-valued rational function $F(s)$, see
\eqref{StateSpace} and let $L_{\rm sq}$ be a realization of the
corresponding squared system $F_{\rm sq}(s)$.
\vskip 0.2cm

For $A$ in \eqref{StateSpace} or in \eqref{eq:HatL}, for
some $q\in[1, n]$ denote
\[
{\rm spect}(A)=\{\lambda_1~,~\ldots~,~\lambda_q\}
\]
with $\lambda_j$ distinct.
\vskip 0.2cm

Consider the following statements

\begin{itemize}
\item[(i)~~~]{} For any $D\in\C^{p\times p}$
\begin{equation}\label{EasySpectASpectL}
{\rm spect}\left(
\begin{smallmatrix}A&B\\ C&D\end{smallmatrix}
\right)\bigcap{\rm spect}(A)=\emptyset .
\end{equation}
\item[(ii)~~]{} There exists a $D\in\C^{p\times p}$ so that
\begin{equation}\label{eq:ExistsD}
{\rm spect}\left(
\begin{smallmatrix}A&B\\ C&D\end{smallmatrix}
\right)\bigcap{\rm spect}(A)=\emptyset .
\end{equation}
\item[(iii)~]{} Consider the realization $L$ of $F(s)$
in \eqref{StateSpace} with $m=p$ (if $m\not=p$ take the squared
counterpart $L_{\rm sq}$ and $F_{\rm sq}(s)$).

\noindent
For each $j=1,~\ldots~,~q$ there
exists $D_j\in\C^{p\times p}$ so that
\[
\lambda_j\not\in{\rm spect}\left(\begin{smallmatrix}
A&B\\ C&D_j\end{smallmatrix}\right).
\]

\item[(iv)~]{} The realization is minimal.
\end{itemize}
Then,
\[
{\rm (i)}~\Longrightarrow~{\rm (ii)}~\Longrightarrow~{\rm (iii)}
~\Longrightarrow~{\rm (iv)}
\]
If $F(s)$ is a scalar rational function, namely $~m=p=1$, then
${\rm (iv)}~\Longrightarrow~{\rm (i)}$.
\end{proposition}

Indeed, the fact that ${\rm (i)}~\Longrightarrow~{\rm (ii)}~\Longrightarrow~{\rm (iii)}$
is straightforward. The equivalence of ${\rm (iii)}$ and ${\rm (iv)}$ is
established in Theorem \ref{Th:Main}.
\vskip 0.2cm

The fact that for scalar systems ${\rm (iv)}~\Longrightarrow~{\rm (i)}$
was first proved (in an elaborate way) in \cite[Theorem 4.1]{MC}. We
next illustrate a straightforward way of showing that.
\vskip 0.2cm

Without loss of generality one can take $A$ to be in its
Jordan canonical form. For example take $n=8$, $q=4$,
\[
L={\footnotesize
\left(\begin{array}{cccccccc|c}
\lambda_1&1        &         &         &         &         &         &         &b_1\\
         &\lambda_1&1        &         &         &         &         &         &b_2\\
         &         &\lambda_1&         &         &         &         &         &b_3\\
         &         &         &\lambda_2&1        &         &         &         &b_4\\
         &         &         &         &\lambda_2&         &         &         &b_5\\
         &         &         &         &         &\lambda_3&1        &         &b_6\\
         &         &         &         &         &         &\lambda_3&         &b_7\\
         &         &         &         &         &         &         &\lambda_4&b_8\\
\hline
c_1      &c_2      &c_3      &c_4      &c_5      &c_6      &c_7      &c_8      &D  
\end{array}\right)},
\]
where zeros were omitted.
\vskip 0.2cm

Using the PBH tests, controllability implies, see e.g. \eqref{eq:PBHcont},
that $\lambda_1$, $\lambda_2$, $\lambda_3$, $\lambda_4$ are
distinct and $b_3$, $b_5$, $b_7$, $b_8$ are non-zero. Observability
implies, see e.g. \eqref{eq:PBHobs}, that $\lambda_1$, $\lambda_2$,
$\lambda_3$, $\lambda_4$ are distinct and $c_1$, $c_4$, $c_6$, $c_8$
are non-zero. Namely, minimality of realization means that
\[
L={\footnotesize
\left(\begin{array}{cccccccc|c}
\lambda_1&1        &         &         &         &         &         &         & * \\
         &\lambda_1&1        &         &         &         &         &         & * \\
         &         &\lambda_1&         &         &         &         &         &\bullet\\
         &         &         &\lambda_2&1        &         &         &         & * \\
         &         &         &         &\lambda_2&         &         &         &\bullet\\
         &         &         &         &         &\lambda_3&1        &         & * \\
         &         &         &         &         &         &\lambda_3&         &\bullet\\
         &         &         &         &         &         &         &\lambda_4&\bullet\\
\hline
\bullet  &  *      &  *      &\bullet  &  *      &\bullet  &  *      &\bullet  &*  
\end{array}\right)},
\]
where $\bullet$ stands for a non-zero scalar and $*$ for ``don't care".
However, looking at $L$ not as a partitioned array, but as a $9\times 9$
matrix, this is exactly the condition for having for $j=1, 2, 3, 4$ the
matrices $(L-\lambda_j{I}_9)$ of a full rank.
(Else, there are rows or columns of zeros). To sum-up, in the context of
scalar systems minimality of realization and condition
\eqref{EasySpectASpectL} are equivalent.
\vskip 0.2cm

The following example illustrates the gap between these conditions
for ~{\em matrix-valued}~ rational functions. Specifically, we
construct realization matrices $A,~B,~C~$ so that in spite of
minimality, for ``many" $D$'s not only the relation \eqref{eq:ExistsD}
does not holds, but in fact \mbox{${\rm spect}(A)\subset{\rm
spect}\left(\begin{smallmatrix}A&B\\ C&D\end{smallmatrix}\right)$}.

\begin{example}\label{Ex:forallD}
{\rm 
Consider the scalar rational functions
$f_1(s)=\frac{2+s}{s}$ and $f_2(s)=\frac{s}{s-2}$
The corresponding minimal realizations are
\[
L_1={\footnotesize\left(\begin{array}{c|c}
0&2\\ \hline 1&1 \end{array}\right)}
\quad\quad\quad\quad
L_2={\footnotesize\left(\begin{array}{c|c}
2&2\\ \hline 1&0 \end{array}\right)}.
\]
Clearly, \mbox{$0={\rm spect}(A_1)\not\in{\rm spect}(L_1)$} and
\mbox{$2={\rm spect}(A_2)\not\in{\rm spect}(L_2)$}.
\vskip 0.2cm

From the above $f_1(s)$ and $f_2(s)$ we now construct the following
matrix valued rational function,
\[
F(s)=\begin{pmatrix}f_1(s)&d_2\\ d_3&f_2(s)\end{pmatrix}=
\begin{pmatrix}\frac{2+s}{s}&d_2\\ d_3&\frac{s}{s-2}\end{pmatrix},
\]
where $d_2, d_3$ are parameters. Its minimal realization is given by
\[
L={\footnotesize\left(\begin{array}{cc|cc}
0&0&2&0\\ 0&2&0&2 \\ 
\hline
1&0&1&d_2\\ 0&1&d_3&1
\end{array}\right)}.
\]
Thus, not only \eqref{EasySpectASpectL} is no longer true, but
in fact for all $d_2, d_3$.
\[
\{ 0, 2\}={\rm spect}(A)\subset
{\rm spect}(L)=\left\{0, ~2, ~1+\sqrt{4+d_2d_3}, ~1-\sqrt{4+d_2d_3}
\right\}.
\]
This example will be further discussed in part II of Example \ref{ExInv}.
}
\qed
\end{example}
\vskip 0.2cm

We conclude this section by pointing out that we do not know whether
or not for $m=p\geq 2$ minimal realization implies that in
\eqref{eq:HatL}
\[
{\rm spect}(A)\bigcap{\rm spect}(\hat{L})
=\emptyset.
\]
\section{
Truncated square realization}
\label{sec:Background}
\setcounter{equation}{0}

In this section we prove Proposition \ref{L:TruncSq}. To this
end we resort to a matrix theory result, whose
interest goes beyond the scope of this work: In
Lemma \ref{L:echelon} below we show that if all we know
about a given matrix is a list of
subsets of its rows which are linearly independent, we
can still impose certain restrictions on its structure.
This is illustrated through a specific, yet rich, example.
\vskip 0.2cm

Let $B$ be a $17\times m~$ matrix, with $~m~$ parameter, be
partitioned to
\begin{equation}\label{eq:partition}
B=\left(\begin{smallmatrix}
B_a\\~\\ \hline\\ B_b\\~\\ \hline\\ B_c
\end{smallmatrix}\right)
\end{equation}
with $B_a\in\C^{6\times m}$, $B_b\in\C^{4\times m}$ and
$B_c\in\C^{7\times m}$. All that is known is that the following
subsets of rows are linearly independent
\begin{equation}\label{eq:LinInd}
(2, 4, 6)\quad{\rm in}\quad B_a~, \quad\quad\quad
(3, 4)\quad{\rm in}\quad B_b~, \quad\quad\quad
(4, 5, 6, 7)\quad{\rm in}\quad B_c~.
\end{equation}
Namely, in $B$ the rows $(2, 4, 6), (9, 10), (14, 15, 16, 17)$
are linearly independent. This implies that $m\geq 4$ and that
${\rm rank}(B)$ can be arbitrary within the range
$[4, \min(m, 17)]$. Indeed, if all rows are spanned by the last
four, numbered $(14, 15, 16, 17)$ then ${\rm rank}(B)=4$. In
contrast $B$ may be of a full rank.
\vskip 0.2cm

Next, one can always find a block diagonal unitary matrix $U$ so that
\[
B=U\tilde{B}
\]
with $\tilde{B}$ in a block echelon form. Namely, a unitary matrix
\[
U={\rm diag}\{U_a\quad U_b\quad U_c\}
\]
with
$U_a\in\C^{6\times 6}$, $U_b\in\C^{4\times 4}$,
$U_c\in\C^{7\times 7}$ and
\begin{equation}\label{eq:TildeB}
\tilde{B}=\left(\begin{smallmatrix}
\tilde{B}_a\\~\\ \hline\\ \tilde{B}_b\\~\\ \hline\\
\tilde{B}_c\end{smallmatrix}\right)
=\left(\begin{smallmatrix}
*&\cdots&\cdots&\cdots&\cdots&\cdots&\cdots&\cdots&\cdots&\cdots&\cdots&*\\
0&\cdots&\cdots&0&\bullet&*&\cdots&\cdots&\cdots&\cdots&\cdots&*\\
*&\cdots&\cdots&\cdots&\cdots&\cdots&\cdots&\cdots&\cdots&\cdots&\cdots&*\\
0&\cdots&\cdots&\cdots&0&\bullet&*&\cdots&\cdots&\cdots&\cdots&*\\
*&\cdots&\cdots&\cdots&\cdots&\cdots&\cdots&\cdots&\cdots&\cdots&\cdots&*\\
0&\cdots&\cdots&\cdots&\cdots&\cdots&0&\bullet&*&\cdots&\cdots&*\\~\\
\hline\\~\\
*&\cdots&\cdots&\cdots&\cdots&\cdots&\cdots&\cdots&\cdots&\cdots&\cdots&*\\
*&\cdots&\cdots&\cdots&\cdots&\cdots&\cdots&\cdots&\cdots&\cdots&\cdots&*\\
0&\cdots&\cdots&\cdots&0&\bullet&*&\cdots&\cdots&\cdots&\cdots&*\\
0&\cdots&\cdots&\cdots&\cdots&\cdots&\cdots&0&\bullet&*&\cdots&*\\~\\
\hline\\~\\
*&\cdots&\cdots&\cdots&\cdots&\cdots&\cdots&\cdots&\cdots&\cdots&\cdots&*\\
*&\cdots&\cdots&\cdots&\cdots&\cdots&\cdots&\cdots&\cdots&\cdots&\cdots&*\\
*&\cdots&\cdots&\cdots&\cdots&\cdots&\cdots&\cdots&\cdots&\cdots&\cdots&*\\
\bullet&*&\cdots&\cdots&\cdots&\cdots&\cdots&\cdots&\cdots&\cdots&\cdots&*\\
0&\cdots&0&\bullet&*&\cdots&\cdots&\cdots&\cdots&\cdots&\cdots&*\\
0&\cdots&\cdots&\cdots&0&\bullet&*&\cdots&\cdots&\cdots&\cdots&*\\
0&\cdots&\cdots&\cdots&\cdots&\cdots&\cdots&0&\bullet&*&\cdots&*
\end{smallmatrix}
\right)
\end{equation}
with $\bullet$ standing for non-zero entry and
$~*$ for ``don't care".
\vskip 0.2cm

Furthermore, for the above matrix $B$, it is always possible to construct 
a full rank~ \mbox{$T_b\in\C^{m\times 4}$} so that, although
$~BT_b~$ has only 4 columns, the list of its linearly
independent rows is as in \eqref{eq:LinInd}, e.g.
\begin{equation}\label{eq:TB}
\tilde{B}T_b=
\left(\begin{smallmatrix}
\tilde{B}_aT_b\\~\\ \hline\\ \tilde{B}_bT_b\\~\\ \hline\\
\tilde{B}_cT_b\end{smallmatrix}\right)
=\left(\begin{smallmatrix}
*&*&*&*\\
0&\bullet&*&*\\
*&*&*&*\\
0&0&\bullet&*\\
*&*&*&*\\
0&0&0&\bullet\\~\\
\hline\\~\\
*&*&*&*\\
*&*&*&*\\
0&0&\bullet&*\\
0&0&0&\bullet\\~\\
\hline\\~\\
*&*&*&*\\
*&*&*&*\\
*&*&*&*\\
\bullet&*&*&*\\
0&\bullet&*&*\\
0&0&\bullet&*\\
0&0&0&\bullet
\end
{smallmatrix}\right).
\end{equation}
More formally, we have the following.

\begin{lemma}\label{L:echelon}
(i) Given $B\in\C^{n\times n}$ nonsingular. One can 
always construct an $n\times n$ unitary matrix $U$ so that 
$UB$ is upper triangular (with non-zero diagonal), i.e.
\[
UB=\hat{B}=\left(\begin{smallmatrix}
\bullet&*&\cdots&\cdots&*\\
0&\bullet&*&\cdots&*\\
0&0&\bullet&*&*\\
\vdots&\vdots&\vdots&\vdots&\vdots\\~\\
0&\cdots&0&\bullet&*\\
0&\cdots&\cdots&0&\bullet
\end{smallmatrix}\right),
\]
with $\bullet$ standing for a non-zero entry and $~*$ for ``don't care".
\vskip 0.2cm

(ii) Given $B\in\C^{n\times m}$ and let $J$ be a subset of
$\{1, 2,~\ldots~,~n\}$. If the rows of $B$ with index $J$ are linearly
independent, then one can always construct an $n\times n$ unitary
matrix $U$ so that in the product matrix $UB$ the rows with index $J$
are in an upper triangular echelon form.
\vskip 0.2cm

(iii) Given
\[
B=\left(\begin{smallmatrix}
B_a\\~\\ \hline\\ B_b\\~\\ \hline\\ B_c\\~\\ \hline\\ \vdots
\end{smallmatrix}\right)
\]
with $B_a\in\C^{n_a\times m}$, $B_b\in\C^{n_b\times m}$,
$B_c\in\C^{n_c\times m}$, $\ldots$ Let $J_a$, $J_b$, $J_c$
$\ldots$ be subsets of $\{1, 2,~\ldots~,~n_a\}$,
$\{1, 2,~\ldots~,~n_b\}$, $\{1, 2,~\ldots~,~n_c\}$
$\ldots$, respectively.
\vskip 0.2cm

If the rows of $B_a$, $B_b$, $B_c$, $\ldots$ with index
$J_a$, $J_b$, $J_c$, $\ldots$, respectively, are linearly independent,
then one can always construct 
unitary matrices 
\mbox{$U_a\in\C^{n_a\times n_a}$}
\mbox{$U_b\in\C^{n_b\times n_b}$}
\mbox{$U_c\in\C^{n_c\times n_c}$} $\ldots$
so that the product matrix
\[
{\rm diag}\{U_a\quad U_b\quad U_c~\ldots~\}B
\]
is in block upper triangular echelon form. Specifically,
in the first $n_a$ rows those with the $J_a$ index are in echelon
form, in the second $n_b$ rows those with the $J_b$ index are in echelon
form, etc.
\vskip 0.2cm

Moreover, if we denote by $\rho$ the maximal cardinality among the sets of
indices $J_a,~J_b,~J_c,~\ldots~$ then $\min(n, m)\geq\rho$ where
$~n:=n_a+n_b+n_c+\ldots$ and ${\rm rank}(B)$ can be arbitrary in the range
$[\rho,~\min(n, m)]$.
\vskip 0.2cm

(iv) In the framework of the previous item, there exists a full
rank matrix $T\in\C^{m\times\rho}$ so that in the $n\times\rho$
product matrix $BT$, the rows with the indices $J_a$, $J_b$, $J_c$,
$\ldots~$ are linearly independent.
\end{lemma}

{\bf Proof}\quad (i)\quad
Assume first that $B$ is $n\times n$ non-singular.  We shall find
it convenient to write down $B$ by its columns as
\[
B=\left(b_1\quad b_2\quad b_3\quad\ldots\quad b_n\right),
\]
where $b_1,~\ldots~,~b_n\in\C^n$ are linearly independent.
Consider the following procedure starting from the left,
\[
\left(\begin{smallmatrix}{b_1}^*\\{b_2}^*\\ \vdots\\b_{n-2}^*\\
b_{n-1}^*\end{smallmatrix}\right)w_n=0\quad\quad
\left(\begin{smallmatrix}b_1^*\\b_2^*\\ \vdots\\b_{n-2}^*\\ {w_n}^*
\end{smallmatrix}\right)w_{n-1}=0\quad\quad
\left(\begin{smallmatrix}b_1^*\\ \vdots\\b_{n-3}\\ 
w_{n-1}^*\\ w_n^*
\end{smallmatrix}\right)w_{n-2}=0\quad\quad\cdots\quad\quad
\left(\begin{smallmatrix}w_2^*\\ w_3^*\\ \vdots\\ w_{n-1}^*\\ w_n^*
\end{smallmatrix}\right)w_1=0.
\]
By construction,
\[
\left(\begin{smallmatrix}w_1^*\\w_2^*\\ \vdots\\ w_{\rho-1}^*\\ w_n^*
\end{smallmatrix}\right)B=
\left(\begin{smallmatrix}
\bullet&*&\cdots&\cdots&*\\
0&\bullet&*&\cdots&*\\
0&0&\bullet&*&*\\
\vdots&\vdots&\vdots&\vdots&\vdots\\~\\
0&\cdots&0&\bullet&*\\
0&\cdots&\cdots&0&\bullet
\end{smallmatrix}\right)
\]
with $\bullet$ standing for a non-zero entry and $~*$ for ``don't care".

Without loss of generality, one can now, normalize
\[
u_j=\frac{w_j}{\| w_j\|_2}
\quad\quad\quad\quad j=1,~\ldots~,~n
\]
\vskip 0.2cm 
so that
\[
U_o=\left(\begin{smallmatrix}u_1^*\\u_2^*\\ \vdots\\ u_{\rho-1}^*\\ u_n^*
\end{smallmatrix}\right)
\]
is unitary and still satisfies
\[
U_oB=\left(\begin{smallmatrix}
\bullet&*&\cdots&\cdots&*\\
0&\bullet&*&\cdots&*\\
0&0&\bullet&*&*\\
\vdots&\vdots&\vdots&\vdots&\vdots\\~\\
0&\cdots&0&\bullet&*\\
0&\cdots&\cdots&0&\bullet
\end{smallmatrix}\right).
\]
\vskip 0.2cm 

(ii)\quad Assume now that $B\in\C^{n\times m}$ and
denote $~r:={\rm rank}(B)$ (clearly \mbox{$\min(m,n)\geq r$}).
\vskip 0.2cm 

Let $~u_{r+1}~,~\ldots~,~u_n\in\C^n$ an orthonormal basis
to the null-space of $B^*$, i.e.
\[
B^*u_j=0
\quad\quad\quad\quad
j=r+1~,~\ldots~,~n.
\]
Let $~u_1~,~\ldots~,~u_r\in\C^n$ an orthonormal completion
of this basis to the whole space. Thus,
\[
U_n:=
\left(u_1\quad u_2\quad u_3\quad\ldots\quad u_n\right),
\]
is a unitary $n\times n$ matrix so that
\[
U_n^*B=\left(\begin{smallmatrix}\tilde{B}\\ 0_{(n-r)\times m}
\end{smallmatrix}\right).
\]
Moreover the linearly independent columns of
\mbox{$\tilde{B}\in\C^{r\times m}$} are numbered from
the left. Namely $\tilde{B}$ can be written as
\[
\tilde{B}=
\left(\begin{smallmatrix}\beta_a&b_1&b_2&\beta_b&\beta_c&b_3&\beta_d&
\ldots&\beta_e&b_\rho&\beta_e&\ldots\end{smallmatrix}\right),
\]
where the columns $b_1,~\ldots~,~b_r\in\C^r$ are
linearly independent and whenever exist, the $\beta$ columns depend
linearly on the $b$ columns to their left (e.g. $\beta_a=0$, both
$\beta_b$ and $\beta_c$ depend on $b_1$ and $b_2$, $\beta_d$ depends
on $b_1, b_2$ $b_3$ etc.). 
Using $U_o\in\C^{r\times r}$ as in
item (i) to construct the matrix
\begin{equation}\label{eq:u}
U=\left(\begin{smallmatrix}U_o&0\\ 0&I_{n-\rho}\end{smallmatrix}\right)
U_n^*
\end{equation}
completes the construction. See part (a) of Example \ref{ExB} for
illustration.
\vskip 0.2cm 

(iii)\quad This follows from the block diagonal structure of $U$.
In each diagonal block, the procedure in item (ii) is applied. See
part (b) of Example \ref{ExB} for illustration.
\vskip 0.2cm 

(iv)\quad This is illustrated in part (c) of Example \ref{ExB} below.
\qed
\vskip 0.2cm 

We next illustrate an application of the above lemma.

\begin{example}\label{ExB}
{\rm
{\bf (a)}~~ Let $B\in\C^{6\times m}$, with $m$ parameter, be so that
rows
\[
(2, 4, 6),
\]
are linearly independent.
This implies that $m\geq 3$ and ${\rm rank}(B)$ can be arbitrary
within the range $[3, \min(m, 6)]$. 
\vskip 0.2cm

We now show that one can always find
a unitary matrix $U$ so that
\[
B=U\tilde{B}
\]
with $\tilde{B}$ in an echelon form. Namely,
\begin{equation}\label{eq:TildeBB}
\tilde{B}=
{\footnotesize
\left(\begin{array}{ccccccccccccc}
*&\cdots&\cdots&\cdots&\cdots&\cdots&\cdots&\cdots&\cdots&\cdots&\cdots&*\\
0&\cdots&\cdots&0&\bullet&*&\cdots&\cdots&\cdots&\cdots&\cdots&*\\
*&\cdots&\cdots&\cdots&\cdots&\cdots&\cdots&\cdots&\cdots&\cdots&\cdots&*\\
0&\cdots&\cdots&\cdots&0&\bullet&*&\cdots&\cdots&\cdots&\cdots&*\\
*&\cdots&\cdots&\cdots&\cdots&\cdots&\cdots&\cdots&\cdots&\cdots&\cdots&*\\
0&\cdots&\cdots&\cdots&\cdots&\cdots&0&\bullet&*&\cdots&\cdots&*
\end{array}\right)}
\end{equation}
with $\bullet$ standing for non-zero entry and
$~*$ for ``don't care".
\vskip 0.2cm

Indeed, first take the following permutation matrix,
\[
P=\left(\begin{smallmatrix}
&1&&&&\\
&&&1&&\\
&&&&&1\\
1&&&&&\\
&&1&&&\\
&&&&1&
\end{smallmatrix}\right),
\]
where zeros were omitted.
Thus, the first three rows of $PB$ are linearly independent.
Next, that ${\rm diag}\{\hat{U}\quad I_3\}$ with $\hat{U}\in\C^{3\times 3}$
as in \eqref{eq:u}. Thus,
\[
\tilde{B}=P^*{\rm diag}\{\hat{U}\quad I_3\}PB.
\]
\vskip 0.2cm

{\bf (b)}~ Let $B\in\C^{17\times m}$, with $m$ parameter,
partitioned as in \eqref{eq:partition} and each block the
linearly independent rows are as in \eqref{eq:LinInd}. We now
show how to obtain a unitary $U$ so that $B=U\tilde{B}$ with
$\tilde{B}$ as in \eqref{eq:TildeB}.
\vskip 0.2cm

Indeed, first take the following permutation matrix,
\[
P={\rm diag}\{P_a\quad P_b\quad P_c\}
\]
with
\[
P_a=\left(\begin{smallmatrix}
&1&&&&\\
&&&1&&\\
&&&&&1\\
1&&&&&\\
&&1&&&\\
&&&&1&
\end{smallmatrix}\right)
\quad\quad
P_b=\left(\begin{smallmatrix}
&I_2\\ I_2&
\end{smallmatrix}\right)
\quad\quad
P_c=\left(\begin{smallmatrix}
&I_4\\ I_3&
\end{smallmatrix}\right)
\]
where zeros were omitted.
%
Thus
\[
PB=\left(\begin{smallmatrix}P_aB_a\\~\\
\hline\\ P_bB_b\\~\\ \hline\\P_cB_c\end{smallmatrix}\right),
\]
where the first three rows in $P_aB_a\in\C^{6\times m}$ are linearly
independent, the first two rows in $P_bB_b\in\C^{4\times m}$ are
linearly independent and the first four rows in $P_cB_c\in\C^{7\times m}$
are linearly independent.
\vskip 0.2cm

Following item ~(a), of this example take now,
\mbox{$U_a={\rm diag}\{\hat{U}_a\quad I_3\}$},
\mbox{$U_b={\rm diag}\{\hat{U}_b\quad I_2\}$},
\mbox{$U_c={\rm diag}\{\hat{U}_c\quad I_3\}$}
where the unitary matrices ~\mbox{$\hat{U}_a\in\C^{3\times 3}$},
\mbox{$\hat{U}_b\in\C^{2\times 2}$} and
\mbox{$\hat{U}_c\in\C^{4\times 4}$} are constructed as in
Lemma \ref{L:echelon}. Thus, in \eqref{eq:TildeB}
\[
\tilde{B}=P^*{\rm diag}\{U_a\quad U_b\quad U_c\}PB.
\]
{\bf (c)}\quad Recall that the number of prescribed linearly
independent rows in $B_a$ was 3, in $B_b$ was 2 and in $B_c$
it was 4. Thus, in the notation of the Lemma $\rho=4$ and we
next show that one can always construct for the matrix $B$
from item ~(b), a full rank $~T_b\in\C^{m\times 4}$ so that
its list of linearly independent rows is as in \eqref{eq:LinInd},
although $~BT_b~$ is only $~17\times 4$.
\vskip 0.2cm

Indeed, if we denote by $e_j$ the standard unit vector where
the $j$-th entry is 1 and zeros elsewhere, for all scalars
$\gamma, \delta$ the $m\times 4$ matrix
\[
T_b=\left(\begin{array}{c|c|c|c}
e_1&e_4+\gamma{e_5}&e_6&e_8+\delta{e_9}
\end{array}\right)
\]
is of a full rank. In fact, one can always choose the scalars
$\gamma, \delta$ so that the product matrix $T_bB$ is as in
\eqref{eq:TB}, so the construction is complete.
\vskip 0.2cm

Finally, we point out that $T_b$ is not unique, for example
one can also take
\[
T_b=\left(\begin{array}{c|c|c|c}
e_1+\gamma{e_8}&e_4&e_6&e_5+\delta{e_9}
\end{array}\right)
\]
with the appropriate $\gamma, \delta$.
}
\qed
\end{example}
\vskip 0.2cm

Thus far for matrix theory results. We now address Proposition
\ref{L:TruncSq} and show that if in \eqref{StateSpace} the
realization is minimal, one can find full rank matrices $T_b$
and $T_c$ so that in \eqref{eq:sq} the realization is of the
smallest possible dimensions (which turn to be square) while
preserving the $A$ matrix and minimality of the realization.
On the way, for a given square matrix $A$, we introduce a
parameterization of all matrices $B$ for which the pairs
$A, B$ are controllable.
\vskip 0.2cm

Recall that the ~{\em geometric multiplicity}~ of
$~\lambda\in{\rm spect}(A)~$ is the number of linearly independent
eigenvectors associated with this $\lambda$, see e.g.
\cite[Definition 1.4.3]{HJ1}. For a given matrix
$A\in\C^{n\times n}$, we shall denote by $\alpha(A)$ the largest
geometric multiplicity among its eigenvalues.
\vskip 0.2cm

Recall also that a matrix $A$ is called ~{\em non-derogatory}~ whenever
$\alpha(A)=1$. (This is equivalent to having the characteristic and the 
minimal polynomials equal, see e.g. \cite[Definitions 1.4.4,
3.2.4.1]{HJ1}, \cite[Corollary 4.4.18]{HJ2}). In turn, this
is closely related to the companion form, see e.g.
\cite[Theorem 3.3.15]{HJ1} (controller form in control engineering
terminology, see e.g. \cite[Section 2.3]{Ka}). For a nice treatment
see \cite[Section 2.2.4]{DP}, \cite[Lemma 6.1.1]{LR}.
\vskip 0.2cm

{\bf Proof of Proposition \ref{L:TruncSq}}\quad 
If a matrix $B\in\C^{n\times m}$ is rank deficient, one can always find
a non-singular $m\times m$ matrix $\tilde{T}$ so that
\[
B=(\tilde{B}\quad 0_{n\times(m-\hat{m})})\tilde{T}
\]
so that $\tilde{B}\in\C^{n\times\hat{m}}$ is of a full rank
$(=\hat{m})$. Moreover, for arbitrary $A\in\C^{n\times n}$, the pairs
$A, B$ and $A, \tilde{B}$ span the same controllable subspace. Thus,
with a slight abuse of notation we shall substitute $\tilde{B}$ by $B$ 
and assume hereafter that
\[
n\geq m={\rm rank}(B).
\]
From the definition of $\alpha=\alpha(A)$,
there exists $\lambda\in{\rm spect}(A)$ so that
\[
n-\alpha={\rm rank}(\lambda{I}_n-A).
\]
Hence,
\[
n-\alpha+m\geq{\rm rank}(\lambda{I}_n-A\quad B),
\]
and if $\alpha>m$, the controllability condition
\eqref{eq:PBHcont} can not hold. Thus, \mbox{$m\geq\alpha$}.
\vskip 0.2cm

To show that one can always choose $T_b$ so that equality
holds, i.e. $m=\alpha$, recall that without loss of generality
one can assume that $A$ is in its Jordan canonical form where
\mbox{$A={\rm diag}\{A_1~,~\ldots~,~A_q\}$} and for
$j=1,~\ldots~,~q$ each $A_j$ (with possibly several Jordan blocks)
contains a single eigenvalue $\lambda_j$ (and
$\lambda_i\not=\lambda_j$ whenever $i\not=j$).
\vskip 0.2cm

To avoid cumbersome notation we proceed by addressing a
specific example. Yet, it is rich enough to cover all cases:
Take $q=3$ and $~n=17$ so that,
\[
A={\rm diag}\{
A_a\quad A_b\quad A_c\}
\]
where for $J_k(\lambda)$ is a $k$-dimensional Jordan block corresponding
to an eigenvalue $\lambda$,
\[
A_a={\rm diag}\{
J_2(\lambda_a),~J_2(\lambda_a),~J_2(\lambda_a)\}
\quad\quad
A_b={\rm diag}\{J_3(\lambda_b),~\lambda_b\}
\quad\quad
A_c={\rm diag}\{J_4(\lambda_c),~\lambda_c{I}_3\}.
\]
%
Then the PBH tests say that controllability of the pair $A, B$ is
equivalent to having $B\in\C^{17\times m}$ ($m$ parameter) whose
subsets of rows indicated in \eqref{eq:LinInd}, are linearly
independent. (As before, this implies that $m\geq 4$ and
${\rm rank}(B)$ can be arbitrary within the range $[4, \min(m, 17)]$). 
\vskip 0.2cm

All such $B$'s are described in part ~(b)~ of Example \ref{ExB}.
In part ~(c)~ of Example \ref{ExB}, a sought $~T_b~$ is constructed.
\vskip 0.2cm

Finally note that the resulting $BT_b$ is of a full rank so indeed
$m=\alpha$.
\vskip 0.2cm

(ii) As observability of the pair $A, C$ is equivalent to controllability
of the pair $A^*, C^*$, this part is omitted. See the illustration below.
\vskip 0.2cm

(iii) Taking $T_b,~T_c$ from items (i), (ii) respectively
completes the construction.
\qed
\vskip 0.2cm

As an illustration we point out that all $C$'s for which the pairs
$A, C$ are observable, where $A\in\C^{17\times 17}$ is as in part
(i) of the proof above, are parameterized as follows:
\vskip 0.2cm

Each $C$ is $p\times 17$ with $p$ parameter and is partitioned to 
\[ 
C=(C_a\quad C_b\quad C_c),
\]
where $C_a\in\C^{p\times 6}$, $C_b\in\C^{p\times 4}$ and
$C_c\in\C^{p\times 7}$ so that the following subsets of columns are
linearly independent
\[
(1, 3, 5)\quad{\rm in}\quad C_a~, \quad\quad\quad
(1, 4)\quad{\rm in}\quad C_b~, \quad\quad\quad
(1, 5, 6, 7)\quad{\rm in}\quad C_c~.
\]
This implies that $p\geq 4$ and that
${\rm rank}(C)$ can be arbitrary within the range $[4, \min(p, 17)]$.
\vskip 0.2cm

We conclude this section by simple illustration of item (iii) in
Proposition \ref{L:TruncSq}

\begin{example}\label{ExNonSq}
{\rm
Given a $1\times 2$-valued rational function (two inputs one
output in control terminology)
\[
F_o(s)=\left(\frac{1}{s}+d~,~\frac{b}{s+1}\right)
\]
where $b, d$ are parameters. It is realized by
\[
L_o=
{\footnotesize
\left(\begin{array}{cr|cc}
0&~~0 &1&0\\
0&-1&0&b\\
\hline
1&~~1 &d&0
\end{array}\right)}.
\]
Take now $T_b=\left(\begin{smallmatrix}1\\1\end{smallmatrix}\right)$
and thus,
\[
F_{\rm sq}(s):=F_o(s)T_b=\frac{1}{s}+\frac{b}{s+1}+d.
\]
The corresponding realization matrix is,
\[
L_{\rm sq}={\footnotesize\left(\begin{array}{cr|c}
0&~~0 &1\\
0&-1&b\\
\hline
1&~~1 &d\end{array}\right)}.
\]
Minimality of both realizations $L_o$ and $L_{\rm sq}$ is equivalent
to $b\not=0$.
\vskip 0.2cm

This example will be further discussed in part I of Example \ref{ExInv}.
}
\qed
\end{example}

\section{proof of Theorem \ref{Th:Main}}
\label{sec:Proof}
\setcounter{equation}{0}

(i)$~\Longrightarrow~$(ii)\\
Recall that we denote by $\alpha$ the largest geometric
multiplicity among the eigenvalues of $A$.
From Lemma \ref{L:TruncSq} it follows that there exist full
rank matrices \mbox{$T_b\in\C^{m\times\alpha}$} and
\mbox{$T_c\in\C^{\alpha\times p}$}
so that $\hat{B}:=BT_b\in\C^{n\times\alpha}$ and
$\hat{C}:=T_cC\in\C^{\alpha\times n}$ are of full rank and
if the realization triple $A, B, C$ was minimal so is
$A, \hat{B}, \hat{C}$.
\vskip 0.2cm

Take now in \eqref{eq:AclA} $K=\eta{T_c}T_b$ with the above
$T_b$, $T_c$ and $\eta>0$ is a scalar parameter. By construction
$K\in\C^{m\times p}$ is of rank $\alpha$ and
\begin{equation}\label{eq:AclNom}
A_{\rm cl}=A+BKC=A+{\eta}\hat{B}\hat{C}.
\end{equation}
We now show that if ${\rm spect}(A)=\{\lambda_1~,~\ldots~,~\lambda_q\}$
by taking $\eta$ ``sufficiently large", $\lambda_j\not\in{\rm spect}(A_{\rm cl})$
for $j=1,~\ldots~,~q$.
\vskip 0.2cm

Let now partition an arbitrary $v_r\in\C^n$ (the subscript stands for ``right")
to \mbox{$v_r=\hat{v}_r+\hat{\hat{v}}_r$}
where $~C\hat{\hat{v}}_r=0$ and $\hat{v}_r$ is in the orthogonal
complement of the null-space of $C$.
\vskip 0.2cm

Note now for a non-zero $\hat{\hat{v}}_r$,
observability implies that
\[
(A_{\rm cl}-\lambda_j{I}_n)\hat{\hat{v}}_r=(A-
\lambda_j{I}_n)\hat{\hat{v}}_r\not=0
\quad\quad j=1,~\ldots~,~q.
\]
On the other hand, for arbitrary non-zero $\hat{v}_r$,
$C\hat{v}_r\not=0$ and by construction 
\mbox{$\hat{B}\hat{C}\hat{v}_r\not=0$}.
Hence for $j=1,~\ldots~,~q$,
\[
\begin{matrix}
\| (A_{\rm cl}-\lambda_jI_n)\hat{v}_r\|&=&\|
(A-\lambda_j{I}_n+{\eta}\hat{B}\hat{C})\hat{v}_r\|
&~\\~
&=&\| (A-\lambda_j{I}_n)\hat{v}_r+
{\eta}\hat{B}\hat{C}\hat{v}_r\|&~\\~
&\geq&|\quad\| (A-\lambda_j{I}_n)\hat{v}_r\|-
{\eta}\|\hat{B}\hat{C}\hat{v}_r\|&
|_{|_{\eta~{\rm sufficiently~large}}}>0.
\end{matrix}
\]
Namely, $(A_{\rm cl}-\lambda_j{I}_n)\hat{v}_r\not=0$.
\vskip 0.2cm

Similarly, partition an arbitrary $v_l\in\C^n$ (the subscript stands for
``left") to \mbox{$v_l=\hat{v}_l+\hat{\hat{v}}_l$}
where ~$\hat{\hat{v}}_l^*B=0$ and $\hat{v}_l$ is in the orthogonal
complement of the null-space of $B^*$.
\vskip 0.2cm

For a non-zero $\hat{\hat{v}}_l$, 
controllability implies that
\[
\hat{\hat{v}}_l^*(A_{\rm cl}-\lambda_jI_n)=
\hat{\hat{v}}_l^*(A-\lambda_jI_n)\not=0
\quad\quad j=1,~\ldots~,~q.
\]
On the other hand, for arbitrary non-zero $\hat{v}_l$,
$\hat{v}_l^*B\not=0$ and by construction
\mbox{$v_l^*\hat{B}\hat{C}\not=0$}.
Hence for $j=1,~\ldots~,~q$,
\[
\begin{matrix}
\|\hat{v}_l^*(A_{\rm cl}-\lambda_jI_n)\|&=&\|
\hat{v}_l^*(A-\lambda_jI_n+{\eta}\hat{B}\hat{C})\|
&~\\~
&=&\|\hat{v}_l^*(A-\lambda_jI_n)+
{\eta}\hat{v}_l^*\hat{B}\hat{C}\|&~\\~
&\geq&|\quad\|\hat{v}_l^*(A-\lambda_jI_n)\|-
{\eta}\|\hat{v}_l^*\hat{B}\hat{C}\|&
|_{|_{\eta~{\rm sufficiently~large}}}>0.
\end{matrix}
\]
Namely, $\hat{v}_l^*(A_{\rm cl}-\lambda_jI_n)\not=0$ so this part of the
claim is established.
\vskip 0.2cm

(ii)$~\Longrightarrow~$(iii)\\
First one can always write
\[
\begin{matrix}
L_{\rm sq}-\lambda{I}_{n+p}&\quad\quad\quad\quad =&
\left(\begin{smallmatrix}A-\lambda{I}_n&~B\\
C&~D-\lambda{I}_p\end{smallmatrix}\right)\\~\\~
&_{|_{\lambda\not\in{\rm spect}(D)}}=&
\left(\begin{smallmatrix}I_n~~&
B(D-\lambda{I}_p)^{-1}\\0&I_p\end{smallmatrix}\right)
\left(\begin{smallmatrix}A_{\rm cl}-\lambda{I_n}&0\\
C&~~~D-\lambda{I}_p\end{smallmatrix}\right),
\end{matrix}
\]
where
\begin{equation}\label{Acl}
A_{\rm cl}:=A+BKC\quad\quad{\rm with}
\quad\quad K:=(\lambda{I}_p-D)^{-1}.
\end{equation}
Thus, the non-singularity of the (different dimensions) matrices
$L_{\rm sq}-\lambda{I}_{n+p}$ and $A_{\rm cl}-\lambda{I_n}$ is
equivalent. Namely, $\lambda\in{\rm spect}(A_{\rm cl})$ if and
only if $\lambda\in{\rm spect}(L_{\rm sq})$.
\vskip 0.2cm

For $j=1,~\ldots~,~q$ take now $D_j=(\lambda_j-\epsilon)I_p$
with $\lambda_j\in{\rm spect}(A)$ and $\epsilon>0$ sufficiently
small. From \eqref{Acl} it follows that that $A_{\rm cl}$ is
of the form \eqref{eq:AclNom} with $\eta={\epsilon}^{-1}$
thus $\lambda_j\not\in{\rm spect}(A_{\rm cl})$ and by the above
construction $\lambda_j\not\in{\rm spect}(L_{\rm sq})$, so 
this part of the claim is established.
\vskip 0.2cm

(iii) $~\Longrightarrow~$ (iv)\quad Trivial.
\vskip 0.2cm

(iv)$~\Longrightarrow~$(i)\\
We find it more convenient to show that if the realization is not
minimal, then
\begin{equation}\label{LambdaInSpectL}
\exists\lambda\in{\rm spect}(A)\quad\quad{\rm so~that}\quad\quad
\lambda\in\bigcap\limits_{D\in\C^{p\times p}}
{\rm spect}\left(\begin{smallmatrix}A&B\\ C&D\end{smallmatrix}\right).
\end{equation}
If the realization $L$ is not controllable, from condition
\eqref{eq:PBHcont} it follows that there exists $\lambda\in\C$
so that
\begin{equation}\label{uncont}
n-1\geq{\rm rank}(\lambda{I}_n-A\quad B).
\end{equation}
This implies that for that same $\lambda$
\begin{equation}\label{SingL}
n+p-1\geq{\rm rank}\left(\begin{smallmatrix}
\lambda{I}_n-A&~~B\\ C&~~*\end{smallmatrix}\right)
\end{equation}
where $*$ stands for ``don't care". Namely, this $\lambda$ is in
${\rm spect}\left(\begin{smallmatrix}A&~~B\\ C&~~*
\end{smallmatrix}\right)$. 
Note now that \eqref{uncont} implies
that, this $\lambda$ is in ${\rm spect}(A)$ (else
${\rm rank}(\lambda{I}_n-A)=n$). Thus, \eqref{LambdaInSpectL} holds.
\vskip 0.2cm

Similarly, if the realization $L$ is not observable, from condition
\eqref{eq:PBHobs} it follows that there exists $\lambda\in\C$
so that
\begin{equation}\label{unobs}
n-1\geq{\rm rank}\left(\begin{smallmatrix}\lambda{I}_n-A\\ C
\end{smallmatrix}\right).
\end{equation}
As before, \eqref{SingL} holds and thus this $\lambda$ is in
${\rm spect}\left(\begin{smallmatrix}A&~~B\\ C&~~*
\end{smallmatrix}\right)$. Note now that \eqref{unobs} implies
that, this $\lambda$ is in ${\rm spect}(A)$ (else
${\rm rank}(\lambda{I}_n-A)=n$). Thus \eqref{LambdaInSpectL} holds
and the proof is complete.
\qed

\section{Families of systems}
\label{sec:Families}
\setcounter{equation}{0}

In linear algebra it is natural to discuss families of matrices
sharing common properties, it is less common address families of
systems. However, the description through the realization matrix
$L$ in \eqref{StateSpace} actually suggests that. Before going
into details, we recall that in \cite{AL}, families of realization
matrices $L$ were studied, in a different framework, by the same
authors. 
\vskip 0.2cm

Let a $(m+n)\times(m+n)$ realization matrix $L$ be given. From
the PBH tests it follows that whenever a realization is not
observable, $L$ has an eigenvector of the
form $\left(\begin{smallmatrix}v_r\\ 0_m\end{smallmatrix}\right)$
for some $0\not=v_r\in\C^n$, i.e. $v_r$ belongs to the orthogonal
complement of the observable subspace associated with the pair
$A, C$. Similarly, whenever a realization is not controllable
$L^*$ has an eigenvector of the
form $\left(\begin{smallmatrix}v_l\\ 0_m\end{smallmatrix}\right)$
for some $0\not=v_l\in\C^n$, i.e. $v_l$ belongs to the orthogonal
complement of the controllable subspace associated with the pair
$A, B$.
\vskip 0.2cm

Hence, it is well known that the controllable (observable) subspace
associated with a given
\[
L={\footnotesize\left(\begin{array}{c|c}
A&B\\ \hline C&D\end{array}\right)}
\]
is identical for all
\[
{\footnotesize\left(\begin{array}{c|c}A+cI_n&B\\ \hline C&*\end{array}\right)}
\quad\quad\quad\forall c\in\C,
\]
where $*$ stands for ``don't care".
This property was already used in this work.
\vskip 0.2cm

Note now that using the same reasoning, whenever $L$ is non-singular also
\[
L^{-1}
\]
shares the same controllable (observable) subspace. From systems point of
view this is not so intuitive since there is no apparent connection between
$F(s)$ realized by $L$ and $F_{\rm inv}(s)$, the
rational function realized by $L^{-1}$.
\vskip 0.2cm
 
Proposition \ref{RationaL} goes along the same lines. It is an immediate
consequence of the following classical matrix theory observation, whose
proof is omitted.

\begin{proposition}
For given: A square matrix  $L$ and a (scalar) polynomial
$\psi(s)$, consider the square matrix $\psi(L)$.
\vskip 0.2cm

The right (left) invariant subspaces of $L$ are contained in the
right (left) invariant subspaces of $\psi(L)$.
\end{proposition}

\vskip 0.2cm

We now illustrate the significance of Proposition \ref{RationaL}.

\begin{example}\label{ExInv}
{\rm
\quad{\bf I.}~
In Example \ref{ExNonSq} we showed that the rational functions
\[
F_o(s)=\left(\frac{1}{s}+d~,~\frac{b}{s+1}\right)
\quad\quad{\rm and}\quad\quad
F_{\rm sq}(s)=\frac{1}{s}+\frac{b}{s+1}+d,
\]
(with $b, d$ parameters) were related. The realization of
$F_{\rm sq}(s)$ was given by,
\[
L_{\rm sq}={\footnotesize\left(\begin{array}{cc|c}
0&~~0 &1\\
0&-1&b\\
\hline
1&~~1 &d\end{array}\right)}.
\]
Consider now the inverse matrix,
\[
L_{\rm sq}^{-1}={\footnotesize\left(\begin{array}{cr|c}
-b-d&1&1\\b&-1&0\\ \hline
1&0&0\end{array}\right)}.
\]
This is a realization of
\[
F_{\rm inv}(s)=\frac{s+1}{s^2+s(b+d+1)+d}~.
\]
Thus, $F_o(s)$, $F_{\rm sq}(s)$ and $F_{\rm inv}(s)$ are all related.
The three respective realizations are minimal, if and only if
$b\not=0$.
\vskip 0.2cm

{\bf II.} Recall that in Example \ref{Ex:forallD} we studied
the realization
\begin{equation}\label{eq:ExL}
L={\footnotesize\left(\begin{array}{cc|cc}
0&0&2&0\\ 0&2&0&2 \\ 
\hline
1&0&1&d_2\\ 0&1&d_3&1
\end{array}\right)},
\end{equation}
(with $d_2, d_3$ parameters) of the rational function,
\[
F(s)=\begin{pmatrix}f_1(s)&d_2\\ d_3&f_2(s)\end{pmatrix}=
\begin{pmatrix}\frac{2+s}{s}&d_2\\ d_3&\frac{s}{s-2}\end{pmatrix}.
\]
We showed that in spite of the minimality of $L$, not only
\mbox{${\rm spect}(A)\bigcap{\rm spect}(L)\not=\emptyset$,}
but in fact \mbox{${\rm spect}(A)\subset{\rm spect}(L)$,}
for all $d_2, d_3$.
\vskip 0.2cm

Consider now the polynomial $\psi(s)=s^2-2s$ and  using
\eqref{eq:ExL}, the corresponding matrix
\begin{equation}\label{eq:Psi(L)}
\psi(L)=L^2-2L=
{\footnotesize\left(\begin{array}{c|c}
\tilde{A}&\tilde{B}\\ \hline
\tilde{B}&\tilde{D}\end{array}\right)},
\end{equation}
with 
\[
\tilde{A}=2I_2\quad\quad\quad
\tilde{B}=2\tilde{C}
=2\left(\begin{smallmatrix}-1&~d_2\\~~d_3&~1\end{smallmatrix}\right)
\quad\quad\quad\tilde{D}=(d_2d_3+1)I_2~.
\]
As, ${\rm spect}(\tilde{A})=\{2,~2\}~$ and $~{\rm spect}(\psi(L))=\{
0,~0,~d_2d_3+3,~d_2d_3+3\}$, it follows that whenever $d_2d_3\not=-1$,
\[
{\rm spect}(\tilde{A})\bigcap{\rm spect}(\psi(L))=\emptyset.
\]
Namely condition \eqref{eq:ExistsD} holds.\quad
%
In fact, $\psi(L)$ is a minimal realization of
\[
F_{\psi}(s)=(d_2d_3+1)\frac{s}{s-2}I_2
\quad\quad\quad\quad d_2d_3\not=-1.
\]
The special case $d_2d_3=-1$ illustrates the that 
in part
(a) of Proposition \ref{RationaL}, (iii) {\em strictly}~
implies (ii).
(The fact that the implication from (i) to (ii) is strict
was addressed in Example \ref{Ex:forallD}).
}
\qed
\end{example}

We believe and hope that this work (along with \cite{AL}) are
just a stage in the study of families systems viewed through
the corresponding realization matrices $L$.
\vskip 0.2cm

\begin{center}
Acknowledgement
\end{center}

\noindent
The authors wish to thank Prof. S. Ter Horst from
the Mathematics group in the
North West University, South Africa,
for pointing out to them reference \cite{MG}.

\end{document}